\documentclass[12pt]{article}

\usepackage{amsmath}
\usepackage{amsthm}
\usepackage{amssymb}
\usepackage{graphics}
\usepackage{graphicx}
\usepackage{graphpap}
\usepackage[arrow,curve,matrix]{xy}

\newtheorem{theorem}{Theorem}
\newtheorem{lemma}[theorem]{Lemma}
\newtheorem{corrol}[theorem]{Corollary}
\newenvironment{jglist}
{\begin{list}{$\bullet$} {
\setlength{\listparindent}{\parindent}\setlength{\parsep}{0 em} }}
{\end{list}}

\def\defn#1{{\bf #1}}

\def\Real{{\mathbb R}}
\def\Cmpx{{\mathbb C}}

\def\Natn{{\mathbb N}}

\def\inter{{\bigcap}}

\def\lldots{,\!..,}

\def\Poisson(#1,#2){{\left\{#1,#2\right\}}}
\def\Set#1{{\left\{#1\right\}}}

\def\ket(#1){{\left|#1\right>}}
\def\bra(#1){{\left<#1\right|}}
\def\ComBrak#1{[#1]}
\def\NonComBrak#1{\langle#1\rangle}
\def\OF#1{{\Aeps^{#1}\!\FAlg}}

\def\piFPn{{\pi_{\FAlg}}}
\def\piPnPM{{{\psi}}}
\def\piFA{{\psi_\FAlg}}

\def\Aelement#1{{\boldsymbol{#1}}}
\def\Aa{\Aelement{a}}
\def\Ab{\Aelement{b}}
\def\Ac{\Aelement{c}}

\def\Af{\Aelement{f}}
\def\Aw{\Aelement{w}}
\def\Ax{\Aelement{x}}
\def\Ay{\Aelement{y}}
\def\Az{\Aelement{z}}
\def\Af{\Aelement{f}}
\def\Ag{\Aelement{g}}
\def\Ah{\Aelement{h}}

\def\Aeps{\Aelement{\varepsilon}}
\def\Ap{\Aelement{p}}
\def\Aq{\Aelement{q}}

\def\AC{\Aelement{C}}
\def\AF{\Aelement{F}}

\def\ffx_#1{\hat\Ax_{#1}}

\def\BPP(#1:#2){{\Big\{#1\,\big|\,#2\Big\}}}
\def\SPP(#1:#2){{\{#1\,|\,#2\}}}

\def\PAlg{{\cal P}(\man)}
\def\Alg{{\cal A}}

\def\FAlg{{\cal F}}
\def\BAlg{{\cal B}}
\def\GAlg{{\cal G}}
\def\JAlg{{\cal J}}
\def\KAlg{{\cal K}}

\def\FOmega{\boldsymbol{\Omega}}

\def\FAlgNil{{\cal F}_0}
\def\FAlgOmega{{\cal F}^\Omega}
\def\FAlgNilOmega{{\cal F}_0^\Omega}

\def\PAlgn{\Cmpx\ComBrak{x_1\lldots x_n}}

\def\EpsAlgInfG{\Cmpx[[\Aeps]]}
\def\EpsAlgInf#1{\EpsAlgInfG\NonComBrak{#1}}

\def\Alg{{\cal A}}

\def\man{{\cal M}}

\def\Cbrak(#1,#2){{C(#1,#2)}}

\title{The Deformation Quantization of $\Real^n$ via the specifiaction
of the commutators}

\author{Jonathan Gratus
\\
\small
Mathematics Department,
School of Informatics
\\
\small
University of Wales, Bangor
Dean Street, Bangor
Gwynedd, UK, LL57 1UT
\\
\small
{email: j@gratus.net, \hspace{5em}
    www.gratus.net}
}

\begin{document}
\maketitle


\begin{abstract}
In a deformation quantization of $\Real^n$, 
the Jacobi identity is automatically
satisfied. This article poses the contrary question: Given a set of
commutators which satisfies the Jacobi identity, is the resulting
associative algebra a deformation quantization of $\Real^n$? It is
shown that the result is true. However care must taken when stating
precisely
how and in which algebra the Jacobi identity is satisfied.
\end{abstract}

MSC 53D55, 81S10, 81R60

\section{Introduction}
\label{ch_intr}

Deformation quantization of a Poisson manifold $\man$ is usually
defined in terms of a star product. 
This is a product on $C^\infty(\man)[[\Aeps]]$ which is the set 
of the infinite sums of complex valued functions of the form
$\sum_{r=0}^\infty \Aeps^r f_r$ where $f_r\in C^\infty(\man)$.
The star product is given by
\begin{equation}
\begin{aligned}
&\star:C^\infty(\man)[[\Aeps]] \times C^\infty(\man)[[\Aeps]]
\to
C^\infty(\man)[[\Aeps]] \;;
\\
&\left(\sum_{r=0}^\infty \Aeps^r f_r \right)
\star
\left(\sum_{s=0}^\infty \Aeps^s f_s \right)
=
\sum_{r=0}^\infty\sum_{s=0}^\infty\sum_{t=0}^\infty
\Aeps^{r+s+t}\, C_t(f_r,g_s)\,,
\end{aligned}
\label{intr_def_star}
\end{equation}
where
\begin{align}
C_r : C^\infty(\man) \times C^\infty(\man) \to C^\infty(\man)\,,
\label{intr_def_star_C_r}
\end{align}
are differential operators,
\begin{align}
C_0 (f,g) = f\, g
\qquad\text{and}\qquad
C_1(f,g) - C_1(g,f) = \Poisson(f,g)\,.
\end{align}
We require $\star$ to be associative,
$\Af\star(\Ag\star\Ah)=(\Af\star\Ag)\star\Ah$ for $\Af,\Ag,\Ah\in
C^\infty(\man)[[\Aeps]]$.

We know from Fedosov, Kontsevich and others that such a construction 
is alway possible, and the degree to which it is unique. For the case that
$\man=\Real^n$ we know that, for a given Poisson structure, it is
unique up to equivalence.
Excellent reviews of deformation quantization are given in
\cite{Sternheimer1,Gutt1}.

\vspace{1 em}

An alternative, and arguably more intuitive, way of defining a
deformation quantization, is via a quotient of the free algebra
$\EpsAlgInf{\Ax_1\lldots\Ax_n}$. In this article we fix $n$ and we write
\begin{align}
\FAlg = \EpsAlgInf{\Ax_1\lldots\Ax_n}\,.
\end{align}
As a set 
\begin{align}
\FAlg 
=
\Set{\sum_{r=0}^\infty \Aeps^r \Aa_r
\,\bigg|\,
\Aa_r\in\Cmpx\NonComBrak{\Ax_1\lldots\Ax_n}
}\,.
\end{align}
Recall $\Cmpx\NonComBrak{\Ax_1\lldots\Ax_n}$ is the free associative
noncommutative algebra generated by $\Set{\Ax_1\lldots\Ax_n}$,
i.e. $\Aa_r$ is a finite sum of finite strings of $\Ax_i$'s, 
with the product of two string given by concatenation.
The deformation parameter $\Aeps\in\FAlg$
commutes with all the elements $\Af\in\FAlg$.
Thus the product in $\EpsAlgInf{\Ax_1\lldots\Ax_n}$ is given by
\begin{align}
\left(\sum_{r=0}^\infty \Aeps^r \Aa_r\right)
\left(\sum_{s=0}^\infty \Aeps^s \Ab_s\right)
=
\sum_{r=0}^\infty\sum_{s=0}^\infty \Aeps^{r+s}\Aa_r\Ab_s\,.
\end{align}
There is the natural surjective homomorphism
\begin{align}
\piFPn:\FAlg\to\PAlgn\;;\qquad
\piFPn(\Aeps)=0\,,\quad
\piFPn(\Ax_i)=x_i\,,
\end{align}
where
$\PAlgn$ is the algebra of polynomials in
$x_1\lldots x_n$.
Recall $\psi:\FAlg\to\GAlg$ is a homomorphism of algebras over
$\Cmpx$ if
$\psi(\Af+\Ag)=\psi(\Af)+\psi(\Ag)$, 
$\psi(\Af\Ag)=\psi(\Af)\psi(\Ag)$ and
$\psi(\lambda)=\lambda$ for $\lambda\in\Cmpx$.

We use the notation of bold symbols for elements of noncommutative
algebras and unbold symbols for elements of commutative algebras.

\vspace{1 em}

In this article we are concerned with the deformation quantization of
$\man=\Real^n$. However for this introduction we consider the more general
case when $\man$ is an algebraic manifold or variety with an algebraic
Poisson structure.  Let $\man$ be an $m$ dimensional manifold or
variety given by
\begin{equation}
\begin{aligned}
\man=\Set{(x_1\lldots x_n)\in\Real^n \bigg| F_s(x_1\lldots x_n)=0 ,\,
  s=1\lldots n\!-\!m }
\end{aligned}
\label{intr_gen_man}
\end{equation}
where $F_s(x_1\lldots x_n)$ are polynomials. 

Let $\PAlg\subset C^\infty(\man)$ be the
subalgebra of polynomial functions in $x_1\lldots x_n$, together with the
projection
\begin{align}
\piPnPM:\PAlgn\to\PAlg=\PAlgn/\Set{F_s=0}\,.
\end{align}
It is convenient to use the notation `$/\Set{F_s=0}$' for $\piPnPM$,
especially in commutative diagrams. Likewise for the other quotienting
homomorphisms.

Let the Poisson structure
$\Poisson(\bullet,\bullet)$ be given by
$\Poisson(x_i,x_j)=C_{ij}(x_1\lldots x_n)$, where
$C_{ij}(x_1\lldots x_n)$ are also polynomials, which in general may
depend on all the $x_i$'s. To be consistent with (\ref{intr_gen_man})
we require $\Poisson(x_i,F_s)=0$.  

A deformation quantization of $(\man,\Poisson(\bullet,\bullet))$ is
given by a choice of $\AF_s,\AC_{ij}\in\FAlg$ where 
\begin{align}
\piFPn(\AF_s)=F_s\,,\quad
\piFPn(\AC_{ij})=C_{ij}\,,\quad
\text{for }
s=1\lldots n\!-\!m
\text{ and }
i,j=1\lldots n\,,
\end{align}
such that the following
diagram of homomorphism of complex associative algebras commutes:
\begin{align}
\xymatrix @C10pc @R4pc {
\FAlg 
\ar[r]^{\big/\Set{\AF_s = 0,\; [\Ax_i,\Ax_j] = \Aeps\AC_{ij}}}
_{\text{i.e. }\piFA}
\ar[d]_(0.3){\big/\Set{\Aeps=0,\; [\Ax_i,\Ax_j] = 0}}
_(0.5){\Ax_i\mapsto x_i} 
_(0.7){\text{i.e. } \piFPn}
&
\Alg
\ar[d]
^(0.3){\big/\Set{\Aeps=0}}
^(0.5){\Ax_i\mapsto x_i} 
^(0.7){\text{i.e. } \pi}
\\
\PAlgn
\ar[r]^{\big/\Set{F_s = 0}}
_{\text{i.e. }\piPnPM}
&
\PAlg
}
\label{DAlg_com_diag}
\end{align}
where
\begin{align}
\piFA:\FAlg\to\Alg=
\FAlg
\big/\Set{
[\Ax_i,\Ax_j] = \Aeps\AC_{ij}
 \,,\;
\AF_s = 0\,,\;
i,j=1\lldots n\,,\;
s=1\lldots n-m
}
\label{intr_A_rels}
\end{align}
and
\begin{align}
\pi:\Alg\to \PAlg\,;\quad \pi(\Aeps)=0\,\quad\pi(\Ax_i)=x_i\,.
\label{intr_pi}
\end{align}
We impose the property on $\Alg$:
\begin{align}
\text{Given\quad $\Af\in\Alg$\quad such that\quad $\Aeps\Af=0$ \quad
  then\quad $\Af=0$.}
\label{intr_A_eps_f_zero}
\end{align}


We use square
brackets to represent the commutator $[\Af,\Ag]=\Af\Ag-\Ag\Af$. Clearly
$\AC_{ij}=-\AC_{ji}$.
It is easy to see that this gives the Poisson 
structure via
\begin{align}
\Poisson({\pi(\Af)},{\pi(\Ag)})=
\pi\left(\Aeps^{-1}[\Af,\Ag]\right)\,.
\label{intr_poisson_form}
\end{align}
Note that the element $\Aeps^{-1}\notin\Alg$, since otherwise the map
$\pi$ would not exist. What we mean by (\ref{intr_poisson_form}) is
that we manipulate $[\Af,\Ag]$ using (\ref{intr_A_rels}) so that it is
of the form $[\Af,\Ag]=\Aeps\Ah$, then
$\Poisson({\pi(\Af)},{\pi(\Ag)})=\pi(\Ah)$. It is not hard to show
that this manipulation is always possible and that the result is independent
of the choice of $\Af$ and $\Ag$.

Borrowing the language from quantum algebra we call an
\defn{ordering}, any linear map
$\omega:\PAlg\to\Alg$ such that $\pi(\omega(f))=f$ for all
$f\in\PAlg$. Orderings are far from unique. For example for the
deformation quantization of $\Real^2$ we can set
$\omega(x_1x_2)=\Ax_1\Ax_2$ (normal order) or
$\omega(x_1x_2)=\tfrac12(\Ax_1\Ax_2+\Ax_2\Ax_1)$ (Wick order) or
$\omega(x_1x_2)=\Ax_1\Ax_2+\Aeps\Ah$ where $\Ah\in\Alg$ is any element.

Given an ordering, we can use this to construct the star product
(\ref{intr_def_star}) via
\begin{align*}
C_r(f,g)=\pi\left(
\Aeps^{-r}\left(
\omega(f)\omega(g)-
\sum_{s=0}^{r-1}
\Aeps^s\omega(C_s(f,g))\right)\right)\,.
\end{align*}
However in general this star product will not
be a differential star product, i.e the $C_r$ 
will not be differential operators.


We observe that $\pi:\Alg\to \PAlg$ must be surjective. This is
because the maps $\piFPn$ and $\piPnPM$ are both surjective.  This
imposes severe restrictions on the possible choices of $\AF_s$ and
$\AC_{ij}$. In general random choices of $\AF_s$ and $\AC_{ij}$ will
either produce inconsistencies or counter (\ref{intr_A_eps_f_zero}).
For example
$\FAlg=\EpsAlgInf{\Ax,\Ay}\,,\;
[\Ax,\Ay]=\Aeps\AC_{\Ax\Ay}=\Aeps\,,\;
\AF_1=\Ay\Ax=0\,,$
then $0=(\Ay\Ax)\Ay=\Ay(\Ax\Ay)=\Aeps\Ay$ in $\Alg$. Thus
$\Ay=0$ and $\pi$ is not surjective.

Examples of deformation quantization constructed in this manner include
\begin{jglist}
\item
The Heisenberg algebra:
\begin{align*}
\FAlg=\EpsAlgInf{\Ap,\Aq}\,,\
[\Ap,\Aq]=i\Aeps\ .
\end{align*}

\item
The Manin Plain:
\begin{align*}
\FAlg=\EpsAlgInf{\Ax,\Ay}\,,\
[\Ax,\Ay]=i\Aeps(\Ax\Ay+\Ay\Ax)\, .
\end{align*}
Note setting $\Aq=(1-i\Aeps)(1+i\Aeps)^{-1}$ gives $\Ax\Ay=\Aq\Ay\Ax$.

\item
The Fuzzy or Noncommutative Sphere \cite{Gratus1}:\,
\begin{align*}
\FAlg=\EpsAlgInf{\Ax,\Ay,\Az}\,,\
&[\Ax,\Ay]=i\Aeps\Az\,,\
[\Ay,\Az]=i\Aeps\Ax\,,\
[\Az,\Ax]=i\Aeps\Ay\,,\
\Ax^2+\Ay^2+\Az^2=1\,.
\end{align*}

\item
The Noncommutative Sphere-Torus \cite{Gratus2}:\,
\begin{align*}
\FAlg=\EpsAlgInf{\Ax,\Ay,\Az}\,,\
&[\Ax,\Ay] = i\Aeps\Az\,,\
[\Ay,\Az] = i\Aeps(\Aw\Ax + \Ax\Aw)\,,\
[\Az,\Ax] = i\Aeps(\Aw\Ay + \Ay\Aw)\,,\
\\
&\Az^2+\Aw^2 = 1
\end{align*}
where $\Aw = \Ax^2 + \Ay^2 - R$ and $R\in\Real$.

\end{jglist}


From now on, we concern ourselves only with the case that
$\man=\Real^n$, so there are no $\AF_s$'s. Clearly since $\Alg$ is an
associative algebra then the Jacobi identity is satisfied,
\begin{align}
[\Ax_i,[\Ax_j,\Ax_k]]+
[\Ax_j,[\Ax_k,\Ax_i]]+
[\Ax_k,[\Ax_i,\Ax_j]]=0
\end{align}
which implies
\begin{align}
[\Ax_i,\AC_{jk}]+
[\Ax_j,\AC_{ki}]+
[\Ax_k,\AC_{ij}]=0\,.
\label{intr_Jac}
\end{align}
We can ask the reverse question, that is: 
\begin{quote}
Given a choice of $\AC_{ij}$
which obey the Jacobi identity (\ref{intr_Jac}), 
is the resulting quotient algebra
\begin{align*}
\piFA:\FAlg\to\Alg=
\FAlg
\big/\Set{
[\Ax_i,\Ax_j] = \Aeps\AC_{ij}\,,\; 
i,j=1\lldots n
}
\end{align*}
a deformation
quantization? 
\end{quote}
We prove that this is true.

However the statement of the theorem is a little tricky since we must
ask in which algebra the Jacobi identity (\ref{intr_Jac}) is being
evaluated. There is no point evaluating it in $\FAlg$ since this will almost
never be satisfied, even when the resulting $\Alg$ is a deformation
quantization of $\Real^n$. Also we cannot evaluate it in $\Alg$ since
(\ref{intr_Jac}) will always be satisfied even if $\Alg$ is not a
deformation quantization of $\Real^n$.

In order to state the theorem we
first define the normal ordered elements of $\FAlg$ and the complete
ordering map.

We define the \defn{normal ordered}
elements of $\FAlg$ as the subset 
\begin{align}
\FAlgOmega=\bigoplus_{s=0}^{\infty}\Aeps^s \FAlgNilOmega \subset \FAlg
\qquad
\text{where}
\qquad
\FAlgNilOmega=\text{span}
\Set{\Ax_1^{r_1}\Ax_2^{r_2}\cdots \Ax_n^{r_n}
\ \Big|\ r_1\lldots r_n\in\Natn}\,.
\label{DAlg_def_F_0_Omega}
\end{align}
There exists a \defn{complete ordering map},
\begin{align}
\phi_\infty:\FAlg\to\FAlgOmega\,,
\end{align}
which is linear and  has the property
\begin{align}
\phi_\infty([\Ax_i,\Ax_j])=\Aeps\phi_\infty(\AC_{ij})\,.
\end{align}
The term ordering is used because again, in some sense, we are choosing
an order. However this is a map on a different space to $\omega$
defined above.
\begin{align}
\xymatrix @C2pc{
{\FAlgOmega}
\ar@{^{(}->}[r]
&
{\FAlg}
\ar[r]^{\piFA} 
\ar[rd]_{\piFPn}
\ar@/^1pc/[l]^{\phi_\infty}
&
{\Alg}
\ar[d]^{\pi}
\\
&&
{\PAlgn}
\ar@/_1pc/[u]_{\omega}
}
\end{align}
This map is defined in section \ref{ch_DAlg}. Using
this map we can define a product on $\FAlgOmega$ via
\begin{align}
\mu : \FAlgOmega\times\FAlgOmega\to\FAlgOmega\;;
\qquad
\mu(\Af,\Ag)=\phi_\infty(\Af\Ag)
\end{align}
This article is to prove the following theorem:

\begin{theorem}
The following are equivalent:
\begin{align}
\bullet \qquad 
&
\text{\parbox{0.8\textwidth}{
The quotient algebra $\Alg$ is a
deformation quantization of $\PAlgn$.
}}
\label{intr_thm_quot}
\\[1em]
\bullet \qquad 
&
\text{\parbox{0.8\textwidth}{The product $\mu$ on $\FAlgOmega$ is
associative.
}}
\label{intr_thm_mu}
\\[1em]
\bullet \qquad
&
\phi_\infty\Big([\Ax_i,\AC_{jk}]+[\Ax_j,\AC_{ki}]+[\Ax_k,\AC_{ij}]\Big)=0
\,,\qquad
\forall i,j,k=1\lldots n
\,.
\label{intr_thm_Jacobie}
\end{align}
\end{theorem}

The proof of Theorem 1 is given in section \ref{ch_DAlg}.  


A simple Corollary to theorem 1 is
\begin{corrol}
Given any $\AC_{\Ax\Ay}\in\FAlg=\EpsAlgInf{\Ax,\Ay}$, then we have a
deformation quantization of $\Real^2$. 
\end{corrol}

One application of Theorem 1 is the establishment of the deformation
quantization of a general manifold or variety $\man$, 
given $\AF_s$ and $\AC_{ij}$. In general this is a difficult task. 
However we can divide the task in two by first
quotienting by $[\Ax_i,\Ax_j]=\Aeps\AC_{ij}$ and then by $\AF_s=0$ as
follows.
\begin{align}
\xymatrix @C6pc @R4pc {
\FAlg 
\ar[r]^{\big/\Set{[\Ax_i,\Ax_j] = \Aeps\AC_{ij}}}
\ar[dr]^{\piFPn} 
&
\BAlg
\ar[r]^{\big/\Set{\AF_s = 0}}
\ar[d]^(0.4){\big/\Set{\Aeps=0}}
\ar[d]^(0.6){\Ax_i\mapsto x_i}
&
\Alg
\ar[d]^{\pi}
\\
&
\PAlgn
\ar[r]^{\big/\Set{F_s = 0}}
&
\PAlg
}
\label{DAlg_com_diag_two}
\end{align}
The task now reduces to: first establishing that the Jacobi identity
given by $\AC_{ij}$ reduces to zero (\ref{intr_thm_Jacobie}) and thus
showing that $\BAlg$ is a deformation quantization of
$\PAlgn$; second: showing that $\AF_s$ is in the centre $\BAlg$,
i.e. that $[\AF_s,\Ax_i]=0$.


\section{Definition of $\phi_\infty$ and proof of
  (\ref{intr_thm_Jacobie})$\Longrightarrow$(\ref{intr_thm_mu}) in Theorem~1}

\label{ch_DAlg}

Let $\FAlgNil=\Cmpx\NonComBrak{\Ax_1\lldots\Ax_n}\subset\FAlg$ be the
algebra of elements with terms with no $\Aeps$ factors. We can also
set $\OF{}=\Set{\Aeps\Af\,|\,\Af\in\FAlg}$, so that
$\FAlg=\FAlgNil\oplus\OF{}$. 
We define the linear map
\begin{align}
\phi:\FAlg\to\FAlgNilOmega\oplus\OF{}
\end{align}
where $\FAlgNilOmega$ is given by (\ref{DAlg_def_F_0_Omega})
as follows: Since $\phi$ is linear we need only define its effect on
words (monomials). 
If $\Af\in\OF{}$ then let $\phi(\Af)=\Af$. If
$\Af\in\FAlgNil$ is a monomial then
$\Af$ is a permutation of
a completely ordered string
\begin{align*}
\overbrace{\Ax_1\cdots\Ax_1}^{\text{$r_1$ factors}}\,
\overbrace{\Ax_2\cdots\Ax_2}^{\text{$r_2$ factors}}
\cdots
\overbrace{\Ax_n\cdots\Ax_n}^{\text{$r_n$ factors}}\,.
\end{align*}
We know, therefore, that $\Af$ can be written in the form
\begin{align*}
\Af=\Ax_1^{r_1}\Ax_2^{r_2}\cdots \Ax_n^{r_n}
+
\text{ terms containing commutators, i.e. terms of 
the form $\Ag[\Ax_i,\Ax_j]\Ah$}
\end{align*}
and we can fix one such expression for each monomial $\Af$.
We then replace
the terms with commutators in using the relation
\begin{align}
[\Ax_i,\Ax_j]\mapsto\Aeps\AC_{ij}
\label{DAlg_com_C}
\end{align}
thus forcing these terms to be in $\OF{}$.
So far we have not completely defined $\phi$.
For example there are two ways of reordering $\Ax_3\Ax_2\Ax_1$:
\begin{equation}
\begin{aligned}
\Ax_3\Ax_2\Ax_1 =&
\Ax_3\Ax_1\Ax_2+\Ax_3[\Ax_2,\Ax_1] =
\Ax_1\Ax_3\Ax_2+\Ax_3[\Ax_2,\Ax_1]+[\Ax_3,\Ax_1]\Ax_2
\\
=&
\Ax_1\Ax_2\Ax_3+\Ax_1[\Ax_3,\Ax_2]+\Ax_3[\Ax_2,\Ax_1]+[\Ax_3,\Ax_1]\Ax_2
\\
\mapsto&
\Ax_1\Ax_2\Ax_3+\Aeps\Ax_1\AC_{32}+\Aeps\Ax_3\AC_{21}+\Aeps\AC_{31}\Ax_2
\end{aligned}
\label{DAlg_reord_321_1}
\end{equation}
and
\begin{equation}
\begin{aligned}
\Ax_3\Ax_2\Ax_1 =&
\Ax_2\Ax_3\Ax_1 + [\Ax_3,\Ax_2]\Ax_1 =
\Ax_2\Ax_1\Ax_3 + \Ax_2[\Ax_3,\Ax_1] + [\Ax_3,\Ax_2]\Ax_1 
\\
=&
\Ax_1\Ax_2\Ax_3 + [\Ax_2,\Ax_1]\Ax_3 + \Ax_2[\Ax_3,\Ax_1] + [\Ax_3,\Ax_2]\Ax_1
\\
\mapsto&
\Ax_1\Ax_2\Ax_3 + \Aeps\AC_{21}\Ax_3 + \Aeps\Ax_2\AC_{31} + \Aeps\AC_{32}\Ax_1
\,.
\end{aligned}
\label{DAlg_reord_321_2}
\end{equation}
We choose the
following algorithm to make $\phi$ well defined. First move all the
$\Ax_1$'s left starting with the left most, 
then move all the $\Ax_2$'s left, and so on. Thus
$\phi(\Ax_3\Ax_2\Ax_1)$ is given by (\ref{DAlg_reord_321_1}) and not
(\ref{DAlg_reord_321_2}).

To make this more precise we define the following relationship on
monomials. Let $\Af\in\FAlgNil$ be a monomial, we write
$\Ax_i\prec\Af$ if all the factors in $\Af$ have a strictly higher
subscript than $\Ax_i$, and $\Ax_i\preceq\Af$ if all the factors in 
$\Af$ have a higher or equal
subscript than $\Ax_i$,
i.e. if $\Af$ has $m$ factors we can write 
$\Af=\Ax_{\sigma(1)}\cdots\Ax_{\sigma(m)}$, where
$\sigma:\Set{1\lldots m}\to\Set{1\lldots n}$. We write
\begin{equation}
\begin{aligned}
\Ax_i\prec\Ax_{\sigma(1)}\cdots\Ax_{\sigma(m)}
\qquad\text{if}\qquad
i<\sigma(k),\ k=1\lldots m\,,
\\
\Ax_i\preceq\Ax_{\sigma(1)}\cdots\Ax_{\sigma(m)}
\qquad\text{if}\qquad
i\le\sigma(k),\ k=1\lldots m\,.
\end{aligned}
\label{DAlg_def_prec}
\end{equation}
Given any monomial $\Af\in\FAlgNil$ there exists a \defn{lowest factor}
with respect to $\prec$ which is the factor $\Ax_i$ such that
$\Ax_i\preceq\Af$. The \defn{lowest factor decomposition} of $\Af$
means writing $\Af=\Ag\Ax_i\Ah$ where $\Ag,\Ah\in\FAlgNil$ are
monomials such that $\Ax_i\prec\Ag$ and $\Ax_i\preceq\Ah$. Thus the
$\Ax_i$ in $\Af=\Ag\Ax_i\Ah$ represents the left most lowest factor of
$\Af$.

Given the monomial $\Ax_{\sigma(1)}\cdots\Ax_{\sigma(m)}$ such that
$\Ax_i\prec\Ax_{\sigma(1)}\cdots\Ax_{\sigma(m)}$ let
\begin{align}
\Cbrak(\Ax_{\sigma(1)}\cdots\Ax_{\sigma(r-1)},\Ax_i)=
\sum_{r=1}^m 
\Ax_{\sigma(1)}\cdots\Ax_{\sigma(r-1)}
\AC_{\sigma(r)\,i}
\Ax_{\sigma(r+1)}\cdots\Ax_{\sigma(m)}
\label{DAlg_def_phibrak}
\end{align}
and let $\Cbrak(1,\Ax_i)=0$.
We observe that $\Cbrak(\bullet,\Ax_i)$ obeys  a Leibniz rule.
\begin{lemma}
Given the monomials $\Af,\Ag\in\FAlgNil$ and $i\in\Set{1\lldots n}$ such
that $\Ax_i\prec\Af\Ag$ then 
\begin{align}
\Cbrak(\Af\Ag,\Ax_i)=\Cbrak(\Af,\Ax_i)\Ag+\Af\Cbrak(\Ag,\Ax_i)
\,.
\label{DAlg_phibrak_libnez}
\end{align}
\end{lemma}
\begin{proof}
First observe that $\Ax_i\prec\Af\Ag$ if and only if both
$\Ax_i\prec\Af$ and $\Ax_i\prec\Ag$. Let 
$\Af=\Ax_{\sigma(1)}\cdots\Ax_{\sigma(m)}$ and 
$\Ag=\Ax_{\sigma(m+1)}\cdots\Ax_{\sigma(s)}$, then
(\ref{DAlg_phibrak_libnez}) follows from (\ref{DAlg_def_phibrak}).
\end{proof}

We set $\phi(1)=1$ and define $\phi(\Af)$ for monomials $\Af\in\FAlgNil$ 
inductively via the lowest factor
decomposition: 
\begin{align}
\phi(\Ag\Ax_i\Ah)=\Ax_i\phi(\Ag\Ah)+\Aeps\Cbrak(\Ag,\Ax_i)\Ah\,,
\label{DAlg_def_phi_induction}
\end{align}
where $\Ax_i\prec\Ag$ and $\Ax_i\preceq\Ah$. Also let
\begin{align}
\Delta(\Af)=\phi(\Af)-\Af\,.
\label{DAlg_def_Delta}
\end{align}


Let $\JAlg$ be the ideal generated by the Jacobi identity, i.e.
\begin{align}
\JAlg=\textup{span}\Set{
\Af\Big([\Ax_i,\AC_{jk}]+[\Ax_j,\AC_{ki}]+[\Ax_k,\AC_{ij}]\Big)\Ag
\ \Big|\ 
\Af,\Ag\in\FAlg,\ 
i,j,k=1\lldots n
}
\label{DAlg_def_Jac_ideal}
\end{align}
and let $\KAlg$ be the algebra generated by the elements
$\Delta(\Ag)$, i.e.
\begin{align}
\KAlg=\textup{span}\Set{
\Af\Delta(\Ag)\Ah
\ \Big|\
\Af,\Ag,\Ah\in\FAlg
}\,.
\label{DAlg_def_Kset}
\end{align}

\begin{lemma}
\label{DAlg_lm_jac}
Let $\Af\in\FAlgNil$
be a monomial such that $\Ax_i\prec\Af$ and $\Ax_j\prec\Af$ then
\begin{align}
  [\AC_{ij},\Af] 
+ [\Cbrak(\Af,\Ax_i),\Ax_j] 
- [\Cbrak(\Af,\Ax_j),\Ax_i] 
\in\JAlg{+}\KAlg\,.
\label{DAlg_phibrak_jacobi}
\end{align}
\end{lemma}
\begin{proof}

We shall prove this via induction on the number of factors in
$\Af$. Clearly if $\Af=1$ then (\ref{DAlg_phibrak_jacobi}) is
trivial. If $\Af=\Ax_k$ then (\ref{DAlg_phibrak_jacobi}) the Jacobi identity.
Assuming (\ref{DAlg_phibrak_jacobi}) is true for $\Af$,
consider for $\Af\to\Ax_k\Af$ gives us:
\begin{align*}
\lefteqn{   
  [\AC_{ij},\Ax_k\Af] 
+ [\Cbrak(\Ax_k\Af,\Ax_i),\Ax_j] 
- [\Cbrak(\Ax_k\Af,\Ax_j),\Ax_i] 
 }
\qquad\qquad&
\\
=&
  [\AC_{ij},\Ax_k]\Af + \Ax_k[\AC_{ij},\Af] 
+ [\AC_{ki}\Af,\Ax_j] + [\Ax_k\Cbrak(\Af,\Ax_i),\Ax_j] 
- [\AC_{kj}\Af,\Ax_i] - [\Ax_k\Cbrak(\Af,\Ax_j),\Ax_i]
\\
=& 
  [\AC_{ij},\Ax_k]\Af 
+ \Ax_k[\AC_{ij},\Af] 
+ [\AC_{ki},\Ax_j] \Af
+ \AC_{ki} [\Af,\Ax_j] 
+ [\Ax_k,\Ax_j] \Cbrak(\Af,\Ax_i)
+ \Ax_k [\Cbrak(\Af,\Ax_i),\Ax_j] 
\\ &
- [\AC_{kj},\Ax_i] \Af
- \AC_{kj} [\Af,\Ax_i] 
- [\Ax_k,\Ax_i] \Cbrak(\Af,\Ax_j)
- \Ax_k [\Cbrak(\Af,\Ax_j),\Ax_i]
\\
=&
 \Big( [\AC_{ij},\Ax_k] 
     + [\AC_{ki},\Ax_j] 
     + [\AC_{jk},\Ax_i] \Big) \Af
+ \Ax_k \Big( [\AC_{ij},\Af] 
            + [\Cbrak(\Af,\Ax_i),\Ax_j] 
            - [\Cbrak(\Af,\Ax_j),\Ax_i] \Big)
\\
&
+\Big(\AC_{ki} [\Af,\Ax_j] 
- [\Ax_k,\Ax_i] \Cbrak(\Af,\Ax_j)\Big)
+\Big( [\Ax_k,\Ax_j] \Cbrak(\Af,\Ax_i)
- \AC_{kj} [\Af,\Ax_i]\Big) 
\in\JAlg{+}\KAlg\,,
\end{align*}
since clearly the first term is in $\JAlg$, and the second term is in
$\JAlg{+}\KAlg$ due to the induction assumption. The third and
fourth terms are in $\KAlg$ due to the following:
since $\Ax_j\prec\Af$, then 
$\phi(\Af\Ax_j)=\Ax_j\phi(\Af)+\Aeps\Cbrak(\Af,\Ax_j)$ and
$\phi(\Ax_j\Af)=\Ax_j\phi(\Af)$, hence
$\phi([\Af,\Ax_j])=\Aeps\Cbrak(\Af,\Ax_j)$.
Hence
\begin{align*}
\AC_{ki} [\Af,\Ax_j] 
- [\Ax_k,\Ax_i] \Cbrak(\Af,\Ax_j) 
&=
\AC_{ki} [\Af,\Ax_j] 
- \Aeps\AC_{ki}\Cbrak(\Af,\Ax_j)
+ \Aeps\AC_{ki}\Cbrak(\Af,\Ax_j)
- [\Ax_k,\Ax_i] \Cbrak(\Af,\Ax_j) 
\\
&=
- \AC_{ki} \Delta([\Af,\Ax_j])  
+ \Delta([\Ax_k,\Ax_i]) \Cbrak(\Af,\Ax_j) 
\in\KAlg\,.
\end{align*}
\end{proof}


\begin{lemma}
\label{DAlg_lm_Delta_f_x}
Let $\Af\in\FAlgNil$ be a monomial such that $\Ax_i\prec\Af$ then
\begin{align}
\phi(\Delta(\Af)\Ax_i)
\in
\Aeps(\JAlg{+}\KAlg)\,.
\label{DAlg_Delta_f_x}
\end{align}
\end{lemma}

\begin{proof}
Let $\Af$ have $m$ factors. We define $\Af_{m-r}$ and $\Ag_r$,
inductively on $r=0\lldots m$, to be monomials in $\FAlgNil$ with $m-r$
and $r$ factors respectively. Let $\Af_{m}=\Af$ and $\Ag_0=1$.
For each $k$ let $\Aa,\Ab\in\FAlgNil$ so that $\Af_{m-r}=\Aa\Ax_k\Ab$
is the lowest factor decomposition of $\Af_{m-r}$. Let
$\Af_{m-r-1}=\Aa\Ab$ and $\Ag_{r+1}=\Ag\Ax_k$. Thus $\Ag_r$ is normal
ordered. We shall prove by induction (backwards) on $r$ that
\begin{align}
\Ah_r\in\Aeps(\JAlg{+}\KAlg)
\qquad\text{where}\qquad
\Ah_r=\phi(\Ag_r\Delta(\Af_{m-r})\Ax_i)\,.
\label{DAlg_Delta_f_r}
\end{align}
For $r=m$ then $\Af_{0}=1$ and $\phi(1)=1$, so $\Delta(1)=0$ and thus
$\Ah_m=0$. Similarly if $r=m-1$ then $\Af_{1}=\Ax_s$ for some $s$ and
$\phi(\Ax_s)=\Ax_s$, so $\Delta(\Ax_s)=0$, thus $\Ah_{m-1}=0$.

Assume $\Ah_{r+1}\in\Aeps(\JAlg{+}\KAlg)$. Let $\Af_{m-r}=\Aa\Ax_k\Ab$
be the lowest factor decomposition of $\Af_{m-r}$. Then
\begin{align*}
\phi(\Af_{m-r})=\phi(\Aa\Ax_k\Ab)=
\Ax_k\phi(\Aa\Ab)+\Aeps\Cbrak(\Aa,\Ax_k)\Ab=
\Ax_k\phi(\Af_{m-r-1})+\Aeps\Cbrak(\Aa,\Ax_k)\Ab\,,
\end{align*}
thus
\begin{equation}
\begin{aligned}
\phi(\Ag_{r} \Delta(\Af_{m-r}) \Ax_i)
=&
\phi(
\Ag_{r} \phi(\Af_{m-r}) \Ax_i)
-
\phi(\Ag_{r} \Af_{m-r} \Ax_i)
\\
=&
\phi(\Ag_r\Ax_k\phi(\Af_{m-r-1})\Ax_i) 
+
\Aeps\Ag_{r}\Cbrak(\Aa,\Ax_k)\Ab\Ax_i
-
\phi(\Ag_{r}\Af_{m-r}\Ax_i)
\\
=&
\phi(\Ag_{r+1}\Delta(\Af_{m-r-1})\Ax_i) 
+
\phi(\Ag_{r+1}\Af_{m-r-1}\Ax_i) 
+
\Aeps\Ag_{r}\Cbrak(\Aa,\Ax_k)\Ab\Ax_i
-
\phi(\Ag_{r}\Af_{m-r}\Ax_i)
\\
=&
\Ah_{r+1}
+
\phi(\Ag_{r}\Ax_k\Aa\Ab\Ax_i) 
+
\Aeps\Ag_{r}\Cbrak(\Aa,\Ax_k)\Ab\Ax_i
-
\phi(\Ag_{r}\Aa\Ax_k\Ab\Ax_i)\,.
\end{aligned}
\label{DAlg_g_phi_x}
\end{equation}
Looking at the terms in (\ref{DAlg_g_phi_x}) we have
\begin{align*}
\phi(
\Ag_{r}\Ax_k \Aa\Ab \Ax_i)
&=
\Ax_i\phi(
\Ag_{r}\Ax_k \Aa\Ab )
+
\Aeps\Cbrak(\Ag_{r}\Ax_k\Aa\Ab,\Ax_i)
\\
&=
\Ax_i\Ag_{r}\Ax_k 
\phi(\Aa\Ab) 
+
\Aeps(
\Cbrak(\Ag_{r},\Ax_i)\Ax_k\Aa\Ab
+
\Ag_{r}\AC_{ki}\Aa\Ab
+
\Ag_{r}\Ax_k\Cbrak(\Aa,\Ax_i)\Ab
+
\Ag_{r}\Ax_k\Aa\Cbrak(\Ab,\Ax_i)
)
\end{align*}
and
\begin{align*}
\phi(\Ag_{r} \Aa\Ax_k\Ab \Ax_i)
=&
\Ax_i\phi(\Ag_{r}\Aa\Ax_k\Ab)+
\Aeps\Cbrak(\Ag_{r}\Aa\Ax_k\Ab,\Ax_i)
\\
=&
\Ax_i\Ag_{r}\phi(\Aa\Ax_k\Ab)+
\Aeps\Cbrak(\Ag_{r}\Aa\Ax_k\Ab,\Ax_i)
\\
=&
\Ax_i\Ag_{r}\Ax_k\phi(\Aa\Ab)+
\Aeps\Cbrak(\Ag_{r}\Aa\Ax_k\Ab,\Ax_i)+
\Aeps\Ax_i\Ag_{r}\Cbrak(\Aa,\Ax_k)\Ab
\\
=&
\Ax_i\Ag_{r}\Ax_k\phi(\Aa\Ab) 
 + 
\Aeps(
\Cbrak(\Ag_{r},\Ax_i)\Aa\Ax_k\Ab
+
\Ag_{r}\Cbrak(\Aa,\Ax_i)\Ax_k\Ab
+
\Ag_{r}\Aa\AC_{ki}\Ab
\\
&\qquad\qquad\qquad
+
\Ag_{r}\Aa\Ax_k\Cbrak(\Ab,\Ax_i)
+
\Ax_i\Ag_{r}\Cbrak(\Aa,\Ax_k)\Ab )\,,
\end{align*}
thus
\begin{align*}
\lefteqn{
\phi(
\Ag_{r}\Ax_k \Aa\Ab \Ax_i)
+
\Aeps\Ag_{r} \Cbrak(\Aa,\Ax_k)\Ab \Ax_i
-
\phi(\Ag_{r} \Aa\Ax_k\Ab \Ax_i)
}
\qquad\qquad
&
\\
=&
\Aeps\Big(
\Cbrak(\Ag_{r},\Ax_i)\Ax_k\Aa\Ab
+
\Ag_{r}\AC_{ki}\Aa\Ab
+
\Ag_{r}\Ax_k\Cbrak(\Aa,\Ax_i)\Ab
+
\Ag_{r}\Ax_k\Aa\Cbrak(\Ab,\Ax_i)
+
\Ag_{r} \Cbrak(\Aa,\Ax_k)\Ab \Ax_i
\\
&
\quad
-
\Cbrak(\Ag_{r},\Ax_i)\Aa\Ax_k\Ab
-
\Ag_{r}\Cbrak(\Aa,\Ax_i)\Ax_k\Ab
-
\Ag_{r}\Aa\AC_{ki}\Ab
-
\Ag_{r}\Aa\Ax_k\Cbrak(\Ab,\Ax_i)
-
\Ax_i\Ag_{r}\Cbrak(\Aa,\Ax_k)\Ab 
\Big)
\\
=&
\Aeps\Big(
\Cbrak(\Ag_{r},\Ax_i)[\Ax_k,\Aa]\Ab
+
\Ag_{r}[\AC_{ki},\Aa]\Ab
-
\Ag_{r}[\Cbrak(\Aa,\Ax_i),\Ax_k]\Ab
+
\Ag_{r}[\Ax_k,\Aa]\Cbrak(\Ab,\Ax_i)
\\
&
\quad
+
\Ag_{r} \Cbrak(\Aa,\Ax_k)[\Ab,\Ax_i]
+
\Ag_{r}[\Cbrak(\Aa,\Ax_k),\Ax_i]\Ab
-
[\Ax_i,\Ag_{r}]\Cbrak(\Aa,\Ax_k)\Ab 
\Big)
\\
=&
\Aeps
\Ag_{r}\Big([\AC_{ki},\Aa]
-
[\Cbrak(\Aa,\Ax_i),\Ax_k]
+
[\Cbrak(\Aa,\Ax_k),\Ax_i]\Big)\Ab
\\
&
+
\Aeps\Big(
[\Ag_{r},\Ax_i]\Cbrak(\Aa,\Ax_k)\Ab 
-
\Cbrak(\Ag_{r},\Ax_i)[\Aa,\Ax_k]\Ab
+
\Ag_{r} \Cbrak(\Aa,\Ax_k)[\Ab,\Ax_i]
-
\Ag_{r}[\Aa,\Ax_k]\Cbrak(\Ab,\Ax_i)
\Big)\,.
\end{align*}
Since $\Ax_i\prec\Aa$ and $\Ax_k\prec\Aa$ then, from lemma
\ref{DAlg_lm_jac}, we have
\begin{align*}
\Ag_{r}\Big([\AC_{ki},\Aa]
-
[\Cbrak(\Aa,\Ax_i),\Ax_k]
+
[\Cbrak(\Aa,\Ax_k),\Ax_i]\Big)\Ab
\in
\JAlg{+}\KAlg\,.
\end{align*}
Also
\begin{align*}
\lefteqn{
[\Ag_{r},\Ax_i]\Cbrak(\Aa,\Ax_k)\Ab
-
\Cbrak(\Ag_{r},\Ax_i)[\Aa,\Ax_k]\Ab
}\qquad &
\\
&=
[\Ag_{r},\Ax_i]\Cbrak(\Aa,\Ax_k)\Ab
-
\Aeps\Cbrak(\Ag_{r},\Ax_i)\Cbrak(\Aa,\Ax_k)\Ab
-
\Cbrak(\Ag_{r},\Ax_i)[\Aa,\Ax_k]\Ab
+
\Cbrak(\Ag_{r},\Ax_i)\Aeps\Cbrak(\Aa,\Ax_k)\Ab
\\
&=
-\Delta([\Ag_{r},\Ax_i])\Cbrak(\Aa,\Ax_k)\Ab
+\Cbrak(\Ag_{r},\Ax_i)\Delta(\Aa,\Ax_k)\Ab
\in
\KAlg\,.
\end{align*}
Hence we have
\begin{align}
\phi(\Ag_{r}\Ax_k \Aa\Ab \Ax_i)
+
\Aeps\Ag_{r} \Cbrak(\Aa,\Ax_k)\Ab \Ax_i
-\phi(\Ag_{r} \Aa\Ax_k\Ab \Ax_i)
\in
\Aeps(\JAlg{+}\KAlg)\,.
\label{DAlg_h_gxabx}
\end{align}
Substituting (\ref{DAlg_h_gxabx}), and the induction assumption
$\Ah_{r+1}\in\Aeps(\JAlg{+}\KAlg)$, into (\ref{DAlg_g_phi_x}) gives
\begin{align}
\Ah_{r}=
\phi(\Ag_{r} \Delta(\Af_r) \Ax_i)
\in
\Aeps(\JAlg{+}\KAlg)\,.
\end{align}
Hence $\Ah_{r}\in\Aeps(\JAlg{+}\KAlg)$ for all $k$.
In particular $\Ah_0\in\Aeps(\JAlg{+}\KAlg)$ 
which is equivalent to (\ref{DAlg_Delta_f_x}).
\end{proof}
\begin{theorem}
\label{DAlg_thm_phi_K}
\begin{align}
\phi(\KAlg)\subset
\Aeps(\JAlg{+}\KAlg)\,.
\label{DAlg_phi_K}
\end{align}
\end{theorem}
\begin{proof}
We need to show $\phi(\Af\Delta(\Ag)\Ah)\in
\Aeps(\JAlg{+}\KAlg)$. From linearity we can assume that
$\Af,\Ag,\Ah\in\FAlg$ are monomials with no $\Aeps$ factor.
We shall use induction on the total number of factors in $\Af$, $\Ag$,
and $\Ah$. Let $\Ax_k$ be the lowest factor in the monomial
$\Af\Ag\Ah$. Thus the lowest factor decomposition of $\Af\Ag\Ah$ is
either: $\Af\Ag\Ah=\Aa\Ax_k\Ab\Ag\Ah$ or $\Af\Ag\Ah=\Af\Aa\Ax_k\Ab\Ah$
or $\Af\Ag\Ah=\Af\Ag\Aa\Ax_k\Ab$ depending on where the left most
$\Ax_k$ occurs. Taking each case in turn:

If $\Af=\Aa\Ax_k\Ab$ then our induction assumption is 
$\phi(\Aa\Ab\Delta(\Ag)\Ah)\in\Aeps(\JAlg{+}\KAlg)$ and
\begin{align}
\phi(\Af\Delta(\Ag)\Ah) =&
\phi(\Aa\Ax_k\Ab\Delta(\Ag)\Ah) 
=
\Ax_k\phi(\Aa\Ab\Delta(\Ag)\Ah) 
+
\Aeps\Cbrak(\Aa,\Ax_k)\Ab\Delta(\Ag)\Ah
\in\Aeps(\JAlg{+}\KAlg)\,.
\label{DAlg_f_axb}
\end{align}

If $\Ag=\Aa\Ax_k\Ab$ then our induction assumption is 
$\phi(\Af\Delta(\Aa\Ab)\Ah)\in\Aeps(\JAlg{+}\KAlg)$ and so
\begin{equation}
\begin{aligned}
\phi(\Af\Delta(\Ag)\Ah) =&
\phi(\Af\Delta(\Aa\Ax_k\Ab)\Ah) 
=
\phi(\Af\phi(\Aa\Ax_k\Ab)\Ah-\Af\Aa\Ax_k\Ab\Ah)
\\
=&
\phi(\Af\Ax_k\phi(\Aa\Ab)\Ah-\Af\Aa\Ax_k\Ab\Ah)+
\Aeps\Af\Cbrak(\Aa,\Ax_k)\Ab\Ah
\\
=&
\Ax_k\phi(\Af\phi(\Aa\Ab)\Ah-\Af\Aa\Ab\Ah)
+
\Aeps(
\Cbrak(\Af,\Ax_k)\phi(\Aa\Ab)\Ah
-
\Cbrak(\Af,\Aa\Ax_k)\Ab\Ah
+
\Af\Cbrak(\Aa,\Ax_k)\Ab\Ah
)
\\
=&
\Ax_k\phi(\Af\Delta(\Aa\Ab)\Ah)
+
\Aeps\Cbrak(\Af,\Ax_k)\Delta(\Aa\Ab)\Ah
\in\Aeps(\JAlg{+}\KAlg)\,.
\end{aligned}
\label{DAlg_g_axb}
\end{equation}

If $\Ah=\Aa\Ax_k\Ab$ then our induction assumption is
$\phi(\Af\Delta(\Ag)\Aa\Ab)\in\Aeps(\JAlg{+}\KAlg)$. Let
$\phi(\Ag)=\Ag_\Omega+\Aeps\hat\Ag$ where
$\Ag_\Omega\in\FAlgNilOmega$. Then 
\begin{align*}
\phi(\Af\Delta(\Ag)\Ah) =&
\phi(\Af\Delta(\Ag)\Aa\Ax_k\Ab) 
=
\phi(\Af\phi(\Ag)\Aa\Ax_k\Ab-\Af\Ag\Aa\Ax_k\Ab)
\\
=&
\phi(\Af\Ag_\Omega\Aa\Ax_k\Ab
-\Af\Ag\Aa\Ax_k\Ab)
+\Aeps\Af\hat\Ag\Aa\Ax_k\Ab
\\
=&
\Ax_k\phi(\Af\Ag_\Omega\Aa\Ab-\Af\Ag\Aa\Ab)
+ \Aeps\Cbrak(\Af\Ag_\Omega\Aa,\Ax_k) \Ab
- \Aeps\Cbrak(\Af\Ag\Aa,\Ax_k)\Ab
+\Aeps\Af\hat\Ag\Aa\Ax_k\Ab
\\
=& 
\Ax_k\phi(\Af\Delta(\Ag)\Aa\Ab)
- \Aeps\Ax_k\Af\hat\Ag\Aa\Ab
+ \Aeps\Cbrak(\Af\Ag_\Omega\Aa,\Ax_k) \Ab
- \Aeps\Cbrak(\Af\Ag\Aa,\Ax_k)\Ab
+\Aeps\Af\hat\Ag\Aa\Ax_k\Ab
\\
=&
\Ax_k\phi(\Af\Delta(\Ag)\Aa\Ab)
+\Aeps(
\Cbrak(\Af\Ag_\Omega\Aa,\Ax_k)
-
\Cbrak(\Af\Ag\Aa,\Ax_k)
+
[\Af\hat\Ag\Aa,\Ax_k]
)\Ab
\\
=&
\Ax_k\phi(\Af\Delta(\Ag)\Aa\Ab)
+\Aeps\Big(
  \Cbrak(\Af,\Ax_k)\Ag_\Omega\Aa
+ \Af\Cbrak(\Ag_\Omega,\Ax_k)\Aa
+ \Af\Ag_\Omega\Cbrak(\Aa,\Ax_k)
- \Cbrak(\Af,\Ax_k)\Ag\Aa
\\&
- \Af\Cbrak(\Ag,\Ax_k)\Aa
- \Af\Ag\Cbrak(\Aa,\Ax_k)
+[\Af,\Ax_k]\hat\Ag\Aa
+\Af[\hat\Ag,\Ax_k]\Aa
+\Af\hat\Ag[\Aa,\Ax_k]
\Big)\Ab
\\
=&
\Ax_k\phi(\Af\Delta(\Ag)\Aa\Ab)
+\Aeps(\Cbrak(\Af,\Ax_k)(\Ag_\Omega-\Ag)+[\Af,\Ax_k]\hat\Ag)\Aa\Ab
\\&
+\Aeps\Af(\Cbrak(\Ag_\Omega,\Ax_k)-\Cbrak(\Ag,\Ax_k)+[\hat\Ag,\Ax_k])\Aa\Ab
+\Aeps\Af((\Ag_\Omega-\Ag)\Cbrak(\Aa,\Ax_k)+\hat\Ag[\Aa,\Ax_k])\Ab
\,,
\end{align*}
since
\begin{align*}
\Cbrak(\Af,\Ax_k)(\Ag_\Omega-\Ag)+[\Af,\Ax_k]\hat\Ag
=&
\Cbrak(\Af,\Ax_k)(\Ag_\Omega+\Aeps\hat\Ag-\Ag)+
([\Af,\Ax_k]-\Aeps\Cbrak(\Af,\Ax_k))\hat\Ag
\\
=&
\Cbrak(\Af,\Ax_k)\Delta(\Ag)-
\Delta([\Af,\Ax_k])\hat\Ag
\in\KAlg
\end{align*}
and from lemma \ref{DAlg_lm_Delta_f_x}, we have
\begin{align*}
\phi(\Delta(\Ag)\Ax_k) =&
\phi((\Ag_\Omega+\Aeps\hat\Ag-\Ag)\Ax_k) 
=
\Ax_k\phi(\Ag_\Omega)+\Aeps\Cbrak(\Ag_\Omega,\Ax_k)
-
\Ax_k\phi(\Ag)-\Aeps\Cbrak(\Ag,\Ax_k)
+
\Aeps\hat\Ag\Ax_k
\\
=&
\Aeps\Cbrak(\Ag_\Omega,\Ax_k)-\Aeps\Cbrak(\Ag,\Ax_k)
+
\Ax_k\Ag_\Omega
-
\Ax_k(\Ag_\Omega+\Aeps\hat\Ag)
+
\Aeps\hat\Ag\Ax_k
\\
=&
\Aeps\Cbrak(\Ag_\Omega,\Ax_k)-\Aeps\Cbrak(\Ag,\Ax_k)
+
\Aeps[\hat\Ag,\Ax_k]
\end{align*}
and so
\begin{align*}
\Cbrak(\Ag_\Omega,\Ax_k)-\Cbrak(\Ag,\Ax_k)+[\hat\Ag,\Ax_k]
\in(\JAlg{+}\KAlg)\,.
\end{align*}
Therefore
\begin{align}
\phi(\Af\Delta(\Ag)\Ah)=
\phi(\Af\Delta(\Ag)\Aa\Ax_k\Ab)
\in\Aeps(\JAlg{+}\KAlg)\,.
\label{DAlg_h_axb}
\end{align}

Hence from (\ref{DAlg_f_axb},\ref{DAlg_g_axb},\ref{DAlg_h_axb}) we see that 
$\phi(\Af\Delta(\Ag)\Ah)\in\Aeps(\JAlg{+}\KAlg)$ where ever the lowest
factor of $\Af\Ag\Ah$ occurs. 
Equation (\ref{DAlg_phi_K}) follows.
\end{proof}

\subsection*{Definition of $\phi_r$ and $\phi_\infty$}

We are now in a position to define $\phi_r$ and thus $\phi_\infty$.
We define the linear maps 
\begin{align}
\phi_r:\FAlg\to
\bigoplus_{s=0}^{r-1}\Aeps^s\!\FAlgNilOmega\oplus\OF{r}
\end{align}
for $r\in\Natn$ by induction.  
Let $\phi_0(\Af)=\Af$ and let $\phi_1(\Af)=\phi(\Af)$.
For each $r\ge1$ let the components of $\phi_{r}(\Af)$ be given by
\begin{align}
\phi_{r}(\Af)=\Af^\Omega_r+\Aeps^{r}\hat\Af_r
\qquad\text{where}\qquad
\Af^\Omega_r\in\bigoplus_{s=0}^{r-1}\Aeps^s\!\FAlgNilOmega
\text{ and }
\hat\Af_r\in\FAlg\,,
\label{DAlg_phi_r_comp}
\end{align}
i.e. $\Af^\Omega_r$
consists of all the terms in $\phi_{r-1}(\Af)$ with coefficients
$\Aeps^0,\Aeps^1,\ldots\Aeps^{r-1}$. Set
\begin{align}
\phi_{r+1}(\Af)=\Af^\Omega_r+\Aeps^r\phi(\hat\Af_r)\,.
\label{DAlg_phi_r_plus_comp}
\end{align}
This results in $\phi_{r}(\Aeps\Af)=\phi_{r-1}(\Af)$.
The limit of these maps is given by the \defn{complete ordering map}: 
\begin{align}
\phi_\infty:\FAlg\to\FAlgOmega=\bigoplus_{s=0}^{\infty}\Aeps^s
\FAlgNilOmega
\end{align}
where $\phi_\infty(\Af)-\phi_r(\Af)\in\Aeps^r\FAlg$.
Given $\Af\in\FAlg$ then it is easy
to see that
\begin{align}
\phi_\infty(\Af)=0 
\qquad\Longleftrightarrow\qquad 
\phi_r(\Af)\in\OF{r}\quad\forall r\in\Natn\,.
\end{align}
\begin{lemma}
Given $\Af,\Ag,\Ah\in\FAlg$ then 
\begin{align}
\phi_\infty(\Af\phi_\infty(\Ag)\Ah-\Af\Ag\Ah)\in\Aeps\phi_\infty(\JAlg{+}\KAlg)
\,.
\label{DAlg_phiinf_fgh_inf}
\end{align}
\end{lemma}
\begin{proof}
Equation (\ref{DAlg_phiinf_fgh_inf}) is equivalent to showing
\begin{align}
\phi_\infty(\Af\phi_r(\Ag)\Ah-\Af\Ag\Ah)\in\Aeps\phi_\infty(\JAlg{+}\KAlg)
\label{DAlg_phiinf_fgh_fin}
\end{align}
for all $r\in\Natn$. This we show by induction on $r$. From 
theorem \ref{DAlg_thm_phi_K} we have
$\phi(\Af\phi(\Ag)\Ah-\Af\Ag\Ah)\in\Aeps(\JAlg{+}\KAlg)$, 
hence 
$\phi_\infty(\Af\phi_1(\Ag)\Ah-\Af\Ag\Ah)\in\Aeps\phi_\infty(\JAlg{+}\KAlg)$
so (\ref{DAlg_phiinf_fgh_fin}) is true for $r=1$. Assume 
(\ref{DAlg_phiinf_fgh_fin}) is true for $r$, then letting
$\phi_r(\Ag)=\Ag^\Omega_r+\Aeps^{r}\hat\Ag_r$ as in
(\ref{DAlg_phi_r_comp}) then
\begin{align*}
\phi_\infty(\Af\Ag^\Omega_r\Ah-\Af\Ag\Ah)+
\Aeps^r\phi_\infty(\Af\hat\Ag_r\Ah)
\in\Aeps\phi_\infty(\JAlg{+}\KAlg)\,.
\end{align*}
From theorem \ref{DAlg_thm_phi_K} we have
$\phi_\infty(\Af\phi(\hat\Ag_r)\Ah-\Af\hat\Ag_r\Ah)
\in\Aeps\phi_\infty(\JAlg{+}\KAlg)$,
hence
\begin{align*}
\Aeps^r\phi_\infty(\Af\phi(\hat\Ag_r)\Ah-\Af\hat\Ag_r\Ah)\in
\Aeps^{r+1}\phi_\infty(\JAlg{+}\KAlg)
\subset\Aeps\phi_\infty(\JAlg{+}\KAlg)\,.
\end{align*}
Thus
\begin{align*}
\phi_\infty(\Af\phi_{r+1}(\Ag)\Ah-\Af\Ag\Ah)=
\phi_\infty(\Af\Ag^\Omega_r\Ah+\Aeps^r\Af\phi(\hat\Ag_r)\Ah-\Af\Ag\Ah)
\in\Aeps\phi_\infty(\JAlg{+}\KAlg)\,.
\end{align*}
Thus (\ref{DAlg_phiinf_fgh_fin}) and hence (\ref{DAlg_phiinf_fgh_inf}).
\end{proof}

We can now prove the difficult part of theorem 1.
\begin{theorem}
(\ref{intr_thm_Jacobie})$\Longrightarrow$(\ref{intr_thm_mu})
\end{theorem}
\begin{proof}
Given $i,j,k=1\lldots n$ and $\Af,\Ag\in\FAlg$, then from
(\ref{DAlg_phiinf_fgh_inf}), we have
\begin{align*}
\phi_\infty\!
\Big(\Af
\phi_\infty\Big(\![\Ax_i,\AC_{jk}]\!+\![\Ax_j,\AC_{ki}]\!+
\![\Ax_k,\AC_{ij}]\!\Big)\Ag
-
\Af\Big([\Ax_i,\AC_{jk}]\!+\![\Ax_j,\AC_{ki}]\!
+\![\Ax_k,\AC_{ij}]\Big)\Ag\!\Big)\!
\in\Aeps\phi_\infty(\JAlg{+}\KAlg).
\end{align*}
But from (\ref{intr_thm_Jacobie}) we have 
$\phi_\infty\Big([\Ax_i,\AC_{jk}]+[\Ax_j,\AC_{ki}]+[\Ax_k,\AC_{ij}]\Big)=0$,
hence
\begin{align*}
\phi_\infty
\Big(\Af\Big([\Ax_i,\AC_{jk}]+[\Ax_j,\AC_{ki}]+[\Ax_k,\AC_{ij}]\Big)\Ag\Big)
\in\Aeps\phi_\infty(\JAlg{+}\KAlg)
\end{align*}
which, from linearity, implies
\begin{align*}
\phi_\infty(\JAlg)\subset\Aeps\phi_\infty(\JAlg{+}\KAlg)\,. 
\end{align*}
From theorem \ref{DAlg_thm_phi_K} we have
$\phi_\infty(\KAlg)\subset\Aeps\phi_\infty(\JAlg{+}\KAlg)$.
Combining these gives
$\phi_\infty(\JAlg{+}\KAlg)\subset\Aeps\phi_\infty(\JAlg{+}\KAlg)$
which implies
$\phi_\infty(\JAlg{+}\KAlg)=\Set{0}$. In other words given
$\Af,\Ag,\Ah\in\FAlg$ we have
\begin{align*}
\phi_\infty(\Af\phi_\infty(\Ag)\Ah)-\phi_\infty(\Af\Ag\Ah)=0\,.
\end{align*}
This implies
\begin{align*}
\mu(\Af,\mu(\Ag,\Ah))-\mu(\mu(\Af,\Ag)\Ah)
=
\phi_\infty(\Af\phi_\infty(\Ag\Ah))-\phi_\infty(\phi_\infty(\Af\Ag)\Ah)
=
0\,.
\end{align*}
Hence $\mu$ is associative.
\end{proof}

\section{Proof of
  (\ref{intr_thm_quot})$\Longrightarrow$(\ref{intr_thm_Jacobie}) of
  Theorem 1}

Recall that (\ref{intr_thm_quot}) states that $\Alg$ is a deformation
quantization of $\Real^n$.
\begin{lemma}
\label{DAlg_lm_eps_a_eq_0}
Given $\Af\in\Alg$ such that $\pi(\Af)=0$ then there exists
$\Af'\in\Alg$ such that $\Af=\Aeps\Af'$.
\end{lemma}

\begin{proof}
Since $\piFA$ is surjective, there exists $\hat\Af\in\FAlg$ such that
$\piFA(\hat\Af)=\Af$. Since $\pi(\piFA(\hat\Af))=\pi(\Af)=0$ then
$\hat\Af\in\ker(\pi\circ\piFA)$ and so $\hat\Af=\Aeps\Aa+
\sum_{ijk} \Ab^k_{ij} [\Ax_i,\Ax_j] \Ac^k_{ij}$ for some
$\Aa,\Ab^k_{ij},\Ac^k_{ij}\in\FAlg$.
Thus
\begin{align*}
\Af=&
\piFA(\hat\Af)
=
\piFA\Big(\Aeps\Aa+
\sum_{ijk} \Ab^k_{ij} [\Ax_i,\Ax_j] \Ac^k_{ij}
\Big)
=
\Aeps\piFA(\Aa)+
\piFA\Big(\sum_{ijk} \Ab^k_{ij} [\Ax_i,\Ax_j] \Ac^k_{ij}
\Big)
\\
=&
\Aeps\piFA\Big(\Aa+\sum_{ijk} \Ab^k_{ij} \AC_{ij}\Ac^k_{ij}
\Big)\,.
\end{align*}
\end{proof}


Let $\Omega:\PAlgn\to\FAlgNilOmega\subset\FAlgOmega$ be the linear
map given by 
$\Omega(x_1^{r_1}x_2^{r_2}\cdots x_n^{r_n})=
\Ax_1^{r_1}\Ax_2^{r_2}\cdots \Ax_n^{r_n}$. Since
$\pi:\Alg\to\PAlgn$ is a surjective homomorphism, we have
$\pi\circ\piFA\circ\Omega$ is the identity on $\PAlgn$. The following
(non commuting) diagram show these maps:
\begin{align}
\xymatrix @C2pc{
{\FAlgOmega}
\ar@{^{(}->}[r]
&
{\FAlg}
\ar[r]^{\piFA} 
\ar[rd]_{\piFPn}
\ar@/^1pc/[l]^(0.3){\phi_\infty}
&
{\Alg}
\ar[d]^{\pi}
\ar@/_2pc/[ll]_{\FOmega}
\\
&&
{\PAlgn}
\ar@/^2pc/[ull]^{\Omega}
}
\end{align}

Let 
\begin{align}
\FOmega:\Alg\to\FAlgOmega\;;
\qquad
\FOmega(\Af)=\sum_{r=0}^{\infty} \Aeps^r\Af_r
\end{align}
where
\begin{align*}
\Af_r=\Omega\pi\left(\frac{1}{\Aeps^{r}}
\left(\Af-\sum_{s=0}^{r-1}\Aeps^s\piFA(\Af_s)\right)\right)
\in\FAlgNilOmega\,.
\end{align*}
\begin{lemma}
$\FOmega:\Alg\to\FAlgOmega$ is well defined and $\FOmega$ is a left
  inverse of $\piFA$, i.e. $\piFA\circ\FOmega$ is the identity of $\Alg$.
\end{lemma}
\begin{proof}
We show by induction that
\begin{align}
\Af-\sum_{s=0}^{r-1}\Aeps^s\piFA(\Af_s)\in\Aeps^r\Alg\,.
\label{DAlg_f_sum_fr}
\end{align}
It is trivial for $r=0$. Assume (\ref{DAlg_f_sum_fr}) is 
true for $r$ so that $\Af_r\in\FAlgNilOmega$ is defined,
then
\begin{align*}
\pi
\left(\frac{1}{\Aeps^{r}}\left(
\Af-\sum_{s=0}^{r}\Aeps^s\piFA(\Af_s)
\right)\right)
=&
\pi
\left(\frac{1}{\Aeps^{r}}\left(
\Af-\sum_{s=0}^{r-1}\Aeps^s\piFA(\Af_s)
\right)\right)
-
\pi\piFA(\Af_r)
\\
=&
\pi\piFA\Omega\pi
\left(\frac{1}{\Aeps^{r}}\left(
\Af-\sum_{s=0}^{r-1}\Aeps^s\piFA(\Af_s)
\right)\right)
-
\pi\piFA(\Af_r)
\\
=&
\pi\piFA(\Af_r)-\pi\piFA(\Af_r)=0
\end{align*}
hence using lemma \ref{DAlg_lm_eps_a_eq_0}, 
(\ref{DAlg_f_sum_fr}) is true for $r+1$. Thus (\ref{DAlg_f_sum_fr}) is
true and $\FOmega:\Alg\to\FAlgOmega$ is well defined.

Given $\Af\in\Alg$ then $\piFA\FOmega(\Af)=\Af$ since
\begin{align*}
\Af-\piFA\FOmega(\Af)=
\Af-\sum_{s=0}^{r-1}\Aeps^s\piFA(\Af_s)
-\sum_{s=r}^{\infty}\Aeps^s\piFA(\Af_s)\in\Aeps^r\Alg
\end{align*}
for all $r\in\Natn$, and hence $\FOmega$ right inverse of $\piFA$.
\end{proof}

\begin{lemma}
\begin{align}
\FOmega\piFA(\Af)=\phi_\infty(\Af)\,.
\label{DAlg_Om_piFA_eq_phiinf}
\end{align}
\end{lemma}

\begin{proof}
From (\ref{DAlg_def_phibrak}) we have
\begin{align*}
\Aeps\piFA(\Cbrak(\Ax_{\sigma(1)}\cdots\Ax_{\sigma(r-1)},\Ax_i))=&
\sum_{r=1}^m 
\piFA(\Ax_{\sigma(1)}\cdots\Ax_{\sigma(r-1)})
\piFA(\Aeps\AC_{\sigma(r)\,i})
\piFA(\Ax_{\sigma(r+1)}\cdots\Ax_{\sigma(m)})
\\
=&
\sum_{r=1}^m 
\piFA(\Ax_{\sigma(1)}\cdots\Ax_{\sigma(r-1)})
\piFA([\Ax_{\sigma(r)},\Ax_i])
\piFA(\Ax_{\sigma(r+1)}\cdots\Ax_{\sigma(m)})
\\
=&
\piFA([\Ax_{\sigma(1)}\cdots\Ax_{\sigma(r-1)},\Ax_i])
\end{align*}
and so $\piFA([\Af,\Ax_i])=\Aeps\piFA(\Cbrak(\Af,\Ax_i))$ for $\Af\in\FAlg$.

By induction on the number of factors in $\Af\in\FAlg$ we show
\begin{align}
\piFA(\phi(\Af))=\piFA(\Af)\,.
\label{DAlg_piFA_phi}
\end{align}
This is true for $\Af=1$. Let $\Af=\Ag\Ax_i\Ah$ be the lowest factor
decomposition of $\Af$. Then 
\begin{align*}
\piFA(\phi(\Af)) 
=&
\piFA(\phi(\Ag\Ax_i\Ah))
=
\piFA(\Ax_i\phi(\Ag\Ah)) + \piFA(\Aeps\Cbrak(\Ag,\Ax_i)\Ah)
\\
=&
\piFA(\Ax_i)\piFA(\phi(\Ag\Ah)) + \piFA([\Ag,\Ax_i]\Ah)
=
\piFA(\Ax_i\Ag\Ah+[\Ag,\Ax_i]\Ah)
=
\piFA(\Ag\Ax_i\Ah)
=
\piFA(\Af)\,.
\end{align*}
Repeatedly applying (\ref{DAlg_piFA_phi}) gives
$\piFA\phi_\infty(\Af)=\piFA(\Af)$.

Note that the
kernel of $\piFA:\FAlg\to\Alg$ consists of elements of the form
$\Af([\Ax_i,\Ax_j]-\Aeps\AC_{ij})\Ag$. Thus
$\textup{ker}(\piFA)\inter\FAlgOmega=\Set{0}$. Hence if
$\Af\in\FAlgOmega$ then
$\Af-\FOmega\piFA(\Af)\in\textup{ker}(\piFA)\inter\FAlgOmega$ so
$\Af-\FOmega\piFA(\Af)=0$. 

Hence since $\phi_\infty(\Af)\in\FAlgOmega$ we have
\begin{align*}
\phi_\infty(\Af)=\FOmega\piFA\phi_\infty(\Af)=\FOmega\piFA(\Af)\,.
\end{align*}
\end{proof}

\begin{theorem}
(\ref{intr_thm_quot})$\Longrightarrow$(\ref{intr_thm_Jacobie}).
\end{theorem}
\begin{proof}
Given $\Af,\Ag\in\FAlg$ we have 
\begin{align*}
\phi_\infty(\Af\phi_\infty(\Ag))
=&
\FOmega\piFA(\Af\,\FOmega\piFA(\Ag))
=
\FOmega(\piFA(\Af)\,\piFA\FOmega\piFA(\Ag))
=
\FOmega(\piFA(\Af)\,\piFA(\Ag))
\\
=&
\FOmega\piFA(\Af\Ag)
=
\phi_\infty(\Af\Ag)\,.
\end{align*}
Given $i,j,k=1\lldots n$ then
$\phi_\infty([\Ax_j,\Ax_k])=\Aeps\phi_\infty(\AC_{jk})$. Since $\FAlg$
is associative
\begin{align*}
0=&\phi_\infty(0)=
\phi_\infty\Big(
[\Ax_i,[\Ax_j,\Ax_k]]+[\Ax_j,[\Ax_k,\Ax_i]]+[\Ax_k,[\Ax_i,\Ax_j]]\Big)
\\
=&
\phi_\infty\Big(
[\Ax_i,\phi_\infty([\Ax_j,\Ax_k])]+
[\Ax_j,\phi_\infty([\Ax_k,\Ax_i])]+
[\Ax_k,\phi_\infty([\Ax_i,\Ax_j])]\Big)
\\
=&
\Aeps\phi_\infty\Big(
[\Ax_i,\AC_{jk}]+[\Ax_j,\AC_{ki}]+[\Ax_k,\AC_{ij}]\Big)
\end{align*}
and using (\ref{intr_A_eps_f_zero}), (\ref{intr_thm_Jacobie}) follows.
\end{proof}

\section{Proof of
  (\ref{intr_thm_mu})$\Longrightarrow$(\ref{intr_thm_quot}) of 
Theorem 1}

\begin{proof}
If $\mu$ is an associative product, then we can define the map
\begin{align}
\pi:\FAlgOmega\to\PAlgn\;;\qquad
\pi\left(\sum_{r=0}^\infty\Aeps^r\Af_r\right)=\Af_0\,.
\end{align}
This map is clearly surjective. Furthermore if $\mu(\Aeps,\Af)=0$ then 
$0=\phi_\infty(\Aeps\Af)=\Aeps\phi_\infty(\Af)\in\FAlg$ so
$\phi_\infty(\Af)=0$, so (\ref{intr_A_eps_f_zero}) is satisfied.
Thus  $(\FAlgOmega,\mu)$ is a deformation
quantization of $\PAlgn$. 
\end{proof}


\section{Conclusion and Discussion}

We have shown that theorem 1 is true. Thus in order to show that a set
of commutation relations forms an associative algebra it is simply
necessary that the Jacobi identity is satisfied, with the correct form
of ordering the generators. A further task would be to show
that the result is independent of how the generators are commuted.

As stated in the introduction this is the first step to establishing a
deformation quantization $\Alg$ of a general manifold $\man$, as in
diagram (\ref{DAlg_com_diag_two}).
Theorem 1 shows that
$\BAlg$ is a deformation quantization of $\Real^n$. 
The second step is to show that the polynomials $\AF_s$ are in the
centre of the algebra $\BAlg$, i.e. that $[\AF_s,\Ax_i]=0$ for all
$\AF_s$ and $\Ax_i$. This enables us to quotient by $\Set{\AF_s=0}$ to
give $\Alg$.

Since we are commuting generators of a noncommutative algebra in a
prescribed way, this work is related to noncommutative Groebner
bases as in the discussion of orderings. Furthermore the map
$\phi:\FAlg\mapsto\FAlg$ 
is an example of a rewrite system.
For instance, the two methods of commuting $\Ax_3\Ax_2\Ax_1$ to give
$\Ax_1\Ax_2\Ax_3$ (\ref{DAlg_reord_321_1}-\ref{DAlg_reord_321_2}) can
be drawn as a hexagon:
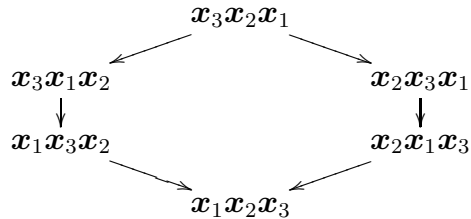
\begin{align*}
\xymatrix @R0.9pc {
&{\Ax_3\Ax_2\Ax_1}\ar[dl]\ar[dr]&
\\
{\Ax_3\Ax_1\Ax_2}\ar[d]&&{\Ax_2\Ax_3\Ax_1}\ar[d]
\\
{\Ax_1\Ax_3\Ax_2}\ar[dr]&&{\Ax_2\Ax_1\Ax_3}\ar[dl]
\\
&{\Ax_1\Ax_2\Ax_3}&
}
\end{align*}
This is an example of a permutahedron. If we calculate the difference
between (\ref{DAlg_reord_321_1}) and (\ref{DAlg_reord_321_2}) we get
the Jacobi identity. So we can say this hexagon represents the
Jacobi identity. If we wish to consider alternative quantization
methods, say working with an associator for non associative algebra,
then we will consider this and more complicated permutahedra.
This has connections to higher dimensional algebras.


\section*{Acknowledgements}
The author would like to thank Arne Sletsj{\o}e, Math Department 
University of Oslo, Tim Porter, Math Department, University of Wales
Bangor, Victor Ufnarovski, Math Department (LTH), Lund University
and the delegates at the Poisson 2004 conference in Luxembourg
for helpful discussions and comments.

The author would like to thank the mathematics departments of the 
University of Oslo and  University of Wales
Bangor for their facilities.

\end{document}